
\baselineskip=14pt
\parskip=10pt
\def\halmos{\hbox{\vrule height0.15cm width0.01cm\vbox{\hrule height
  0.01cm width0.2cm \vskip0.15cm \hrule height 0.01cm width0.2cm}\vrule
  height0.15cm width 0.01cm}}

\magnification=\magstephalf

\def\1{{\overline{1}}}
\def\2{{\overline{2}}}
\parindent=0pt
\overfullrule=0in

\def\frac#1#2{{#1 \over #2}}
\centerline
{\bf 
Experimenting with Ap\'ery Limits and WZ pairs
}
\bigskip
\centerline
{\it Robert DOUGHERTY-BLISS and Doron ZEILBERGER}

\bigskip

\qquad\qquad\qquad {\it  In fond memory of the amazing Borwein brothers: Jonathan (20 May 1951 - 2 August 2016) and Peter (10 May 1953 - 23 August 2020), great pioneers of Experimental Mathematics}

{\bf Abstract}: This article, dedicated with admiration in memory of Jon and Peter Borwein, 
illustrates by example, the power of experimental mathematics, so dear to them both, by experimenting with so-called Ap\'ery limits and WZ pairs.
In particular we prove a weaker form of an intriguing conjecture of Marc Chamberland and Armin Straub (in an article  dedicated to Jon Borwein),
and generate lots of new Ap\'ery limits. We also rediscovered an infinite family of cubic irrationalities, that suggested very good effective irrationality
measures (lower than Liouville's generic $3$), and that we conjectured to go down to the optimal $2$. As it turned out, as pointed
out by Paul Voutier (see the postscript  kindly written by him), our conjectures follow from deep results in number theory.
Nevertheless we believe that further experiments with our Maple programs would lead to new and interesting results.

{\bf Preface: The amazing Borwein Brothers}

[This section was written by DZ, hence the first-person ``I"]

I first met the Borwein brothers in the historic Centenary Ramanujan conference, at the University of Illinois at Urbana-Champaign, that took place June 1-5, 1987. They were already then legendary. They just
published the classic [BoBo1] and their talk [BoBo2] was a natural continuation of their great book, furnishing yet-faster-converging series for computing $\pi$, inspired by the work
of Ramanujan. Here is the conference photo:

{\tt https://sites.math.rutgers.edu/\~{}zeilberg/mamarim/mamarimhtml/Ramanujan100.jpg} \quad \quad .


You can see Jon Borwein in the front row (6th from the right, sporting a Pi T-shirt), and Peter Borwein in the fourth row (7th from the left). Of course I was in awe of them,
and they complemented the other notable $\pi$-brothers, Gregory and David Chudnowsky (who were present, but are not in the photo).

Since then, Jon and Peter became {\it pioneers} in Experimental Mathematics. They founded the {\bf Center for Constructive and Experimental Mathematics} (CECM) at Simon Fraser University

{\tt http://www.cecm.sfu.ca/} \quad ,

that is still flourishing today.

Then Jon, the epitome of {\bf type-A personality}, always with so many plans, went {\it down-under} and founded {\bf CARMA} (COMPUTER-ASSISTED RESEARCH MATHEMATICS AND ITS APPLICATIONS),

{\tt https://carma.newcastle.edu.au/about/} \quad ,

and this center is just as flourishing today.

Together with David Bailey,  Jon wrote the {\bf bible} of Experimental Mathematics [BoBa]. 
This was complemented by a second volume (where  Ronald Girgensohn joined them) [BoBaGi].
These two volumes are so engaging, yet very deep.
See my raving book review for {\it American Scientist} [Ze3]:

{\tt https://sites.math.rutgers.edu/\~{}zeilberg/mamarim/mamarimPDF/mathexp.pdf} \quad .

Note that, defying the almost universal mathematical tradition of naming the authors {\it alphabetically}, Jon Borwein is listed as first author in these two volumes.
I am sure that it is for a very good reason!

I should also mention that both Jon and Peter did so much to cultivate the $\pi$-cult, and the anthology [BerBoBo] is an indispensible source for all $\pi$-lovers.

{\bf The present article: A continuation of Section 10.3 of the Borweins' masterpiece ``Pi and the AGM''}

The present article is inspired by the work of Roger Ap\'ery [Ap]  and Frits Beukers [Beu], beautifully exposited in {\it Pi and the AGM}, Section~10.3. It also continues
the work in [Ze2], and the recent delightful article, by Marc Chamberland and Armin Straub [ChaS], that was dedicated to the memory of Jon Borwein.

{\bf Ap\`ery Limits}

One way that Ap\'ery's seminal proof of the irrationality  of $\zeta(3)$  could have been discovered, in a {\it counterfactual world}, was to
consider, {\it out of the blue}, the second-order recurrence
$$
n^3 u_n \,- \, \left(17 n^{2}+51 n +39\right) \left(2 n +3\right)\,u_{n-1} \,+ \, (n-1)^3 u_{n-2}=0 \quad,
$$
and let $a_n$ and $b_n$ be the solutions of that recurrence with {\bf initial conditions}
$$
a_0=0, a_1=6 \quad; \quad
b_0=1, b_1=5 ,
$$
then let the computer compute many terms, evaluate $\frac{a_{1000}}{b_{1000}}$ to many decimals, and then use Maple's {\tt identify}, and {\it lo and behold}, get that it (most probably) equals $\zeta(3)$
(i.e. $\sum_{i=1}^{\infty} \frac{1}{i^3}$). Then, still {\it empirically} and {\it numerically}, after rewriting $\frac{a_n}{b_n}$ as $\frac{a'_n}{b'_n}$, where now
{\bf both} numerator and denominators are integers (initially $b_n$ are integers, but $a_n$ are not), estimate that there exists a {\it positive} number $\delta$ (about $0.0805$) such that
$$
|\frac{a'_n}{b'_n} - \zeta(3)| \leq \frac{CONSTANT}{(b'_n)^{1+\delta}} \quad ,
$$
that immediately entails (see [vdP]) that $\zeta(3)$ is irrational.

This theme is pursued in [Ze2] and much more recently in [DKZ], where the motivation was to discover irrationality proofs of other constants.

In the above-mentioned delightful article [ChaS], the authors  decided to study Ap\'ery limits independently of their potential for suggesting  irrationality proofs, 
continuing the work of Zudilin [Zu] and Almkvist, van Straten and Zudilin [AlvSZ].
Take {\it any} linear recurrence, that came up `naturally' (e.g. satisfied by a binomial coefficients sum), then define another solution using different
initial conditions, and take the limit of the ratios of the two sequences, and see what happens.

In particular, they made the following intriguing conjecture (Conjecture 9 in [ChaS]).

{\bf Conjecture} (Chamberland and Straub) For $d \geq 3$, the minimal order recurrence satisfied by
$$
A^{(d)}(n) := \sum_{k=0}^{n} \, {{n} \choose {k}}^d \quad,
$$
has a  unique solution, $B^{(d)}(n)$,  with $B^{(d)}(0)=0$ and  $B^{(d)}(1)=1$ that has the property that:
$$
\lim_{ n \rightarrow \infty} \frac{B^{(d)}(n)}{A^{(d)}(n)} \, = \, \frac{\zeta(2)}{d+1} \quad.
$$

Surprisingly, Chamberland and Straub did not address the {\it rate of convergence}. After all, the original motivation of Ap\'ery was to furnish sequences
of rational numbers converging to the desired constant ($\zeta(3)$ or $\zeta(2)$) with {\it exponentially decaying rate}. In fact he needed
more than that (and he lucked out), that the denominators do not grow too fast, but setting this number-theoretical issue aside, even if our sequence of rational numbers
does not imply irrationality, it is desirable to demand that the error in  the approximation decays {\it exponentially}.

It is easy to see that if the Chamberland-Straub conjecture is true, then the following stronger (and far more interesting!) version is also true. Hence
their convergence is equivalent to

{\bf Stronger Chamberland-Straub  Conjecture} For $d \geq 3$, the minimal order recurrence satisfied by
$$
A^{(d)}(n) := \sum_{k=0}^{n} \, {{n} \choose {k}}^d \quad,
$$
has a  unique solution, $B^{(d)}(n)$,  with $B^{(d)}(0)=0$ and  $B^{(d)}(1)=1$ that has the property that there exists  real numbers $C$ and $\alpha>1$ such that
$$
 \lim_{ n \rightarrow \infty} \, |\frac{B^{(d)}(n)}{A^{(d)}(n)} - \frac{\zeta(2)}{d+1}| \leq \frac{C}{\alpha^n} \quad.
$$

Using the results in [Ze1] we can prove a weaker form of this stronger conjecture. 

{\bf Theorem 1}: There exists {\it some} (not necessarily minimal) recurrence and {\it some} initial conditions for the $B^{(d)}(n)$ sequence
such that both $A^{(d)}(n)$ and  $B^{(d)}(n)$  satisfy the {\bf same} linear recurrence equation with polynomial coefficients and
real numbers $C$ and $\alpha>1$ such that
$$
 \lim_{ n \rightarrow \infty} \, | \frac{B^{(d)}(n)}{A^{(d)}(n)} - \frac{\zeta(2)}{d+1}| \leq \frac{C}{\alpha^n} \quad.
$$

In order to prove this, we need to recall the following theorem from [Ze1].
Below $N$ is the forward shift operator, and a general linear recurrence operator, of order $L$, say,  has the form
$$
\sum_{i=0}^{L}\,  c_i(n) \, N^i \quad,
$$
so another way of saying that a sequence $a(n)$ satisfies the linear recurrence equation
$$
\sum_{i=0}^{L}\,  c_i(n) \, a(n+i) \, = \, 0 \quad,
$$
is to say that the sequence $a(n)$ is {\bf annihilated} by the {\bf operator} $\sum_{i=0}^{L}\,  c_i(n) \, N^i$, i.e.
$$
\left ( \sum_{i=0}^{L}\,  c_i(n) \, N^i \right ) a(n) \, = \, 0 \quad.
$$

We need the following theorem.

{\bf Theorem 2} (Theorem 9 of [Ze1]): Let $c(n,k)$ be the potential function of
a WZ 1-form $F(n,k) \delta k + G(n,k) \delta n $ 
in the two variables $(n,k)$. In other words,
 
$$
  F(n,k)= c(n,k+1)-c(n,k) \,, \, \, \,G(n,k)=c(n+1,k)-c(n,k) \, \,, 
$$
 
and let $b(n,k)$ be closed-form (i.e. `proper-hypergeometric'). Let 
$$
  a(n):= \sum_{k=0}^{n} c(n,k) b(n,k)  \,, \, \, \, \, \, \, b(n):= \sum_{k=0}^{n} b(n,k)     \, \,. 
$$
 
 There exist
(rapidly exhibitable) linear recurrence operators
with polynomial coefficients $R (N,n)$ and
$S (N,n)$ such that
 
$$
 R (N,n) b(n)  \,= \,0 \,, \, \, \, \, \, S (N,n) R (N,n) a(n) \,= \,0 \, \, \, \,. 
$$         
 
Furthermore, there exist rapidly exhibitable
closed-form ``certificates" $B(n,k)$ and $D(n,k)$ such that
the following routinely verifiable identities are true:
 
$$
 R (N,n) b(n,k) = B(n,k)-B(n,k-1) \,, 
$$          
 
$$
  S (N,n) R (N,n) ( b(n,k) c(n,k) ) = S (N,n) ( c(n,k)B(n,k) - c(n,k-1)B(n,k-1)) + D(n,k)-D(n,k-1) \, \, \,. 
$$
 
In addition, $B(n,k)/b(n,k)$ and $D(n,k)/(b(n,k)F(n,k))$ are
both rational functions.

{\bf Proof of Theorem 1}:  Consider the potential function, let's call it $c(n,k)$ 
$$
 c(n,k):= 2 \sum_{m=1}^{n} {{(-1)^{m-1}} \over{ m^2 }} + 
\sum_{m=1}^{k} {{(-1)^{n+m-1}} \over { m^2  (  { n \atop m }  )  
          (  { n+m \atop m }  ) }}   \quad,
$$  
that comes from the WZ form
$$
 \omega_{ \zeta (2) }  \,:=  \, {{(-1)^{(n+k)} k!^2 (n-k-1)!}\over{ (n+1)\,(n+k+1)!}}   [ (n+1)\,\delta k  \,+ \, 2(n-k)\, \delta n  \,  ] 
$$        
(see [Ze1]), then it follows from Theorem 2, that for {\it any} {\it binomial coefficients sum}  \quad ,
$$
A(n):=\sum_{k=0}^{n} b(n,k) \quad,
$$
(not just powers of ${{n} \choose {k}}$), defining
$$
B(n)\,=\, \sum_{k=0}^{n} b(n,k)\,c(n,k) \quad ,
$$
that there exists a linear recurrence satisfied by {\it both} sequences $A(n)$ and $B(n)$.

Since $c(n,k)$ converges to $\zeta(2)$ no matter how you approach $(\infty,\infty)$ on the discrete region $\{n \geq k \geq 0\}$,  $B(n)/A(n)$, beings a {\it weighted-average}
of $c(n,k)$ for $0 \leq k \leq n$, also must converge  to $\zeta(2)$.  

It also follows from the proof of Theorem 2 in [Ze1] that the convergence rate is exponentially decaying. \halmos

One of our dreams is to find a {\bf decision procedure} that:

{\bf Inputs}: An {\bf arbitrary} linear recurrence equation with {\bf polynomial coefficients} of order $L$, say, and any two solutions given by {\bf some} initial conditions, 
$$
a(0) \, =\, a_0 \quad, \quad \dots \quad ,\quad a(L-1) \, =\, a_{L-1} \quad,
$$
$$
b(0) \, =\, b_0 \quad,\quad \dots \quad,\quad  b(L-1) \, =\, b_{L-1} \quad ,
$$

{\bf Outputs}: Yes if and only if $\lim_{n \rightarrow \infty} \frac{a(n)}{b(n)} \, = \,0$ \quad .

This decision procedure exists for {\bf constant coefficients} linear recurrences (it goes back at least to Euler), but we have no clue how to do it for linear recurrences with polynomial coefficients.

Assuming that we have such an `{\it oracle}', it would be immediate to prove the Chamberland-Straub conjecture for any specific $d$. Both $A^{(d)}(n)$ and $B^{(d)}(n)$  satisfy the same
{\it minimal recurrence} of $A^{(d)}(n)$. Let  
$$
B'^{(d)}(n)=\sum_{k=0}^{n} {{n} \choose {k}}^r\,c(n,k) \quad ,
$$
then, by {\it linearity},  $B^{(d)}(n) - B'^{(d)}(n)$ also satisfies that very same (non-minimal recurrence), and if we can rigorously decide whether or not
$$
\lim_{ n \rightarrow \infty}  \frac{B^{(d)}(n) - B'^{(d)}(n)}{A^{(d)}(n)} \, = \, 0 \quad ,
$$
we would be able to completely prove the Chamberland-Straub conjecture (for each specific power $d$).

Since we do not (yet!) have such a decision procedure, we can do the next-best thing and prove is empirically. This is done for $3 \leq d \leq 9$ in the output file

{\tt https://sites.math.rutgers.edu/\~{}zeilberg/tokhniot/oAperyWZnew1.txt} \quad,

where the (rigorous!) proof of the weaker form of the Chamberland-Straub conjecture is given, followed by a non-rigorous, `empirical proof' of the full conjecture.

{\bf Having fun with Ap\'ery Limits}

``Ap\'ery's incredible proof appears to be a mixture of miracles and mysteries.
The dominating question is how to generalise all this, down to the Euler
constant $\gamma$ and up to the general $\zeta(t)$?.'' [vdP, sec.~10]

The intervening forty years have not answered van der Poorten's question; we
still do not know how to construct a ``nice'' recurrence to prove the
irrationality of fixed constants. However, as mentioned in [ChaS], the {\it
inverse} problem is just as difficult! That is, given solutions $A(n)$ and
$B(n)$ to some recurrence with polynomial coefficients, if $$\lim_{n
\to \infty} \frac{B(n)}{A(n)}$$ exists, what does it equal? (The same problem
for recurrences with {\it constant} coefficients is much easier.)

The conjecture of [ChaS] considered in the previous section was made with
numerical evidence. In some sense this conjecture was ``easy'' because most
computer algebra systems can identify $r \zeta(k)$ for reasonably simple
rationals $r$ and small $k$. However, for an arbitrary Ap\'ery limit $\lim_{n
\to \infty} B(n) / A(n)$ we are usually stranded.

To help experiment, the Maple package

{\tt https://sites.math.rutgers.edu/\~{}zeilberg/tokhniot/AperyLimits.txt}

systematically searches for recurrences and attempts to guess the corresponding
Ap\'ery limits using our enhanced version of Maple's {\tt identify} command.
The recurrences come from either the Zeilberger algorithm [Ze0] or from the
Almkvist-Zeilberger algorithm [AlZ].

For example, consider the sequence of functions
$$f_n(x) = \left( \frac{(2x + 1)(1 + 3x)}{x} \right)^n x^{-2/3}.$$
The Almkvist-Zeilberger algorithm shows that the sequence $$a(n) = \int_{|x| =1} f_n(x)\ dx$$ satisfies the recurrence

$$a(n) = 9n\frac{(2n - 1)5}{(3n - 1)(3n + 1)} a(n - 1) - 9\frac{n(n-1)}{(3n-1)(3n+1)}a(n
-2).$$

We will now forget $a(n)$ itself and make use only of the recurrence. Let
$A(n)$ and $B(n)$ satisfy this same recurrence with initial conditions
$$A(0) = 1, \quad A(1) = 0; \qquad B(0) = 0, \quad B(1) = 1.$$
Then, empirically,
$$\lim_{n \to \infty} \frac{B(n)}{A(n)} = -\frac{28}{15} 18^{1/3} - \frac{32}{45} 18^{2/3} - \frac{76}{15}.$$

Going through the same process with
$$f_n(x) = \left( \frac{(3x + 1)(1 + 4x)}{x} \right)^n x^{-2/3}$$
gives an empirical limit of
$$-\frac{80}{21} 6^{1/3} - \frac{44}{21} 6^{2/3} - \frac{148}{21}.$$

Both sequences are instances of the more general family
$$f_n(x) = \left( \frac{(cx + 1)(1 + (c + 1)x)}{x} \right)^n x^{-2/3}.$$
Indeed, going through the same steps with arbitrary $c$ suggests that we always
obtain a ``cubic'' Ap\'ery limit. In particular, we conjectured that the induced
Ap\'ery limit is the real root of the cubic
$$
64+144 x \left(1+2 c \right)+108 x^{2} \left(3 c^{2}+3 c +1\right)+27 x^{3} \left(1+2 c \right) = 0,
$$
which we determined thanks to the PSLQ algorithm [BaFe]. 

Even better than the Ap\'ery limit itself, as $c \to \infty$, it seems that the
{\it effective} irrationality measures suggested by our computations go down to
$2$. Recall that for a cubic irrationality, the `vanilla' irrationality measure
promised by Liouville is $3$, so any `infinite family' of cubic irrationalities
that give irrationality measures smaller than $3$ is of interest. Such families
were found by Alan Baker [Bak], Gregory Chudnowsky [Chu],  Michael Bennett [Be], and Paul Voutier [V1][V2].

As  it turned out, our conjectures follow from known deep results in number theory, see the postscript kindly written by Paul Voutier.
But this was just {\it one} example, and hopefully further experiments would lead to new results.

For a computer-generated exploration of this experiment, see the output file

{\tt https://sites.math.rutgers.edu/\~{}zeilberg/tokhniot/oAperyLimits5.txt} \quad .

Similar to this infinite family of cubic irrationalities, the sequence of
functions
$$f_n(x) = \left( \frac{(cx + 1)(1 + (c + 1)x)}{x} \right)^n x^{-1/2}.$$
empirically produces the quadratic Ap\'ery limit
$$
-3c - \frac{3}{2} - 3 \sqrt{c^2+c} \quad .
$$
Every quadratic irrationality has an effective irrationality measure of $2$
(from its continued fraction expansion), but our diophantine approximations are
{\it different} from the standard continued fraction continuants. This is all
conjecture, but we are sure that it won't be too hard to prove it. Since the
stakes are so low, and we are experimental mathematicians, we leave this to the
reader.

For a computer-generated discussion of this experiment, see the file

{\tt https://sites.math.rutgers.edu/\~{}zeilberg/tokhniot/oAperyLimits6.txt} \quad .

Alas, most constants that we were able to identify are known to be irrational
by general theorems (e.g. the  Lindemann-Weierstrass  theorem). For those unidentified constants our results may suggest
explicit irrationality measures, and for the identified constants they may
suggest {\it effective} irrationality measures.

{\bf Postscript by Paul Voutier}

The authors conjecture that the induced Ap\'{e}ry limit is the real root, $\alpha_{c}$,
of the cubic polynomial in $x$
$$
f_{c}(x)=64+144 x \left(1+2 c \right)+108 x^{2} \left(3 c^{2}+3 c +1\right)+27 x^{3} \left(1+2 c \right) = 0
$$
and claim that as $c \rightarrow \infty$, it seems that the effective irrationality
measures suggested by their computations go down to $2$.

We show here that their suggestion is true.
A new such family of algebraic numbers with effective irrationality
measures improving on Liouville's measure would be of considerable
interest. But $\alpha_{c}$ is in
${\bf Q} \left( ((c+1)/c)^{1/3} \right)$ and
effective irrationality measures for elements of such fields have been known
since Baker's pioneering work in 1964 [Bak]. This was later improved by
Chudnowsky [Chu], with explicit versions of Chudnowsky's work produced by Easton
[E], Bennett [Ben] and Voutier [V1].

Using Maple, one can quickly verify that
$$
\frac{4((c+1)/c)^{1/3}-4}{-(9c+3)((c+1)/c)^{1/3}+9c+6}
$$
is a root of $f_{c}(x)$ and, since real when taking the real cube root of
$(c+1)/c$, it must equal $\alpha_{c}$.

Using Theorem~2.1 of [V1] with $a=c+1$ and $b=c$, we can take $d=0$ there,
and so
$$
\kappa=\frac{\log\left( e^{0.911} \left( \sqrt{c+1}+\sqrt{c} \right)^{2} \right)}
{\log\left( e^{-0.911} \left( \sqrt{c+1}-\sqrt{c} \right)^{-2} \right)}.
$$

For $c \geq 4$, we have $\kappa<2$. Thus,
$$
\left| ((c+1)/c)^{1/3} - p/q \right| > \frac{1}{10^{120}(c+1)|q|^{\kappa+1}}
$$
for all integers $p$ and $q$ with $q \neq 0$, which
improves on Liouville's irrationality measure for such $c$.

Since $\alpha_{c}$ is a fractional linear transformation of $((c+1)/c)^{1/3}$
with rational coefficients, we can use Lemma~6.3 of [V2] to obtain an
effective irrationality measure for $\alpha_{c}$ too. In the notation of
Lemma~6.3 of [V2], we have $a_{1}=4$, $a_{2}=-4$, $a_{3}=-(9c+3)$,
$a_{4}=9c+6$ with $\theta=((c+1)/c)^{1/3}$ and $\theta'=\alpha_{c}$. For $c \geq 4$,
it can be shown that
$$
\left| a_{4}+a_{3} \theta \right|
\left( \left| a_{3} \right| \left( 1+\left| \theta' \right| \right) + \left| a_{1} \right| \right)^{\kappa}
< \left| a_{4}+a_{3} \theta \right|
\left( \left| a_{3} \right| \left( 1+\left| \theta' \right| \right) + \left| a_{1} \right| \right)^{2}
< 1000c^{2},
$$
so
$$
\left| \alpha_{c} - p/q \right| > \frac{1}{10^{123}(c+1)^{3}|q|^{\kappa+1}},
$$
for all integers $p$ and $q$ with $q \neq 0$.

Furthermore, as $c \rightarrow \infty$, we have
$$
\kappa \rightarrow 
\frac{\log\left( e^{0.911} 4c \right)}{\log\left( e^{-0.911}c/4 \right)}
\rightarrow 1,
$$
so Dougherty-Bliss and Zeilberger's observation that the effective irrationality
measures approach $2$ as $c \rightarrow \infty$ is true.

{\bf Acknowledgment}: Many thanks are due to Michael Bennett and Wadim Zudilin for very helpful discussions.
Also many thanks to Armin Straub for very helpful comments on an earlier version, and to Paul Voutier who
authored the postscript.

\vfill\eject

{\bf REFERENCES}

[AlZ]  Gert Almkvist and Doron Zeilberger, {\it The method of integrating under the integral sign}, J. Symbolic Comp. {\bf 10} (1990), 571-591.
\hfill\break
{\tt https://sites.math.rutgers.edu/\~{}zeilberg/mamarim/mamarimPDF/duis.pdf} \quad .

[AlvSZ] Gert Almkvist, Duco van Straten, and Wadim Zudilin, {\it Ap\'ery limits of differential equations of orders $4$ and $5$}, Fields Inst. Commun. Ser. {\bf 54} (2008), 105-123.

[Ap] Roger Ap\'ery, 
{\it ``Interpolation de fractions continues et irrationalit\'e
de certaine constantes''}
Bulletin de la section des sciences du C.T.H.S. \#3
p. 37-53, 1981.

[BaFe] Helaman Ferguson and David Bailey, {\it A polynomial time, numerically
stable integer relation algorithm}, RNR Technical Report RNR-91-032.

[Bak] A. Baker,
{\it Rational approximations to $2^{1/3}$ and other algebraic numbers},
Quart. J. Math. Oxford {\bf (2) 15} (1964), 375--383.

[Ben] Michael Bennett, {\it Effective measures of irrationality for certain algebraic numbers},  J. Australian Math. Soc. (series A) {\bf 62} (1997), 329-344.

[Beu] Frits Beukers, {\it A note on the irrationality of $\zeta(2)$ and $\zeta(3)$}, 
Bull. London Math. Soc. {\bf 11} (1979), 268-272.
Reprinted in [BerBoBo], 434-438.

[BerBoBo] Lennard Berggren, Jonathan Borwein, and Peter Borwein, {\it ``Pi: a source book''}, Springer, 1997.

[BoBo1] Jonathan Borwein and Peter Borwein, {\it ``Pi and the AGM''}, Wiley, 1987.

[BoBo2] Jonathan Borwein and Peter Borwein, {\it More Ramanujan-type series for $1/\pi$} in: ``Ramanujan Revisited'', Proceedings of the Centenary Conference, University of Illinois at Urbana-Champaign, June 1-5, 1987,
edited by George E. Andrews et. al., Academic Press.

[BoBa] Jonathan Borwein and David Bailey ,  {\it ``Mathematics by experiment: plausible reasoning in the 21st century''}, second edition, A.K. Peters/CRC Press, 2008.

[BoBaGi] Jonathan Borwein, David Bailey and Roland Girgensohn,  {\it ``Experimentation in mathematics: computational paths to discovery''}, A.K. Peters/CRC Press, 2004.

[ChaS] Marc Chamberland and Armin Straub, {\it Ap\'ery limits: experiments and proofs},   6 Nov 2020. {\tt https://arxiv.org/abs/2011.03400} \quad .

[Chu] G.V. Chudnowsky, {\it On the method of Thue-Siegel}, Ann. of  Math. {\bf 117} (1983), 325-382.

[DKZ] Robert Dougherty-Bliss, Christoph Koutschan, and Doron Zeilberger, {\it Tweaking the Beukers Integrals in search of more miraculous irrationality proofs \`A La Ap\'ery},
submitted. 
\hfill\break
{\tt https://sites.math.rutgers.edu/\~{}zeilberg/mamarim/mamarimhtml/beukers.html} \quad .

[E1] D. Easton,
{\it Effective irrationality measures for certain algebraic numbers},
Math. Comp. {\bf 46} (1986), 613--622.

[vdP] Alf van der Poorten, {\it A Proof that Euler missed...
Ap\'ery's proof of the irrationality of $\zeta(3)$}, Math. Intelligencer {\bf 1} (1979), 195-203.
Reprinted in [BerBoBo], 439-447.

[V1] P. M. Voutier,
{\it Rational approximations to $2^{1/3}$ and other algebraic numbers revisited},
Journal de Th\'{e}orie des Nombres de Bordeaux {\bf 19} (2007), 265--288.

[V2] P. M. Voutier,
{\it Thue's Fundamentaltheorem, I},
Acta Arith. {\bf 143} (2010), 101--144.

[Ze0] Doron Zeilberger, {\it The method of creative telescoping}, 
J. Symbolic Computation {\bf 11} (1991), 195-204.  \hfill\break
{\tt https://sites.math.rutgers.edu/\~{}zeilberg/mamarimY/creativeT.pdf} \quad .

[Ze1] Doron Zeilberger, {\it Closed form (pun intended!)}, 
 in: ``Special volume in memory of Emil Grosswald", M. Knopp
and M. Sheingorn,
eds., Contemporary Mathematics {\bf 143} 579-607, AMS, Providence (1993). \hfill\break
{\tt https://sites.math.rutgers.edu/\~{}zeilberg/mamarim/mamarimhtml/pun.html} \quad ,

[Ze2] Doron Zeilberger, {\it Computerized deconstruction}, Advances in Applied Mathematics {\bf 30} (2003), 633-654. \hfill\break
{\tt https://sites.math.rutgers.edu/\~{}zeilberg/mamarim/mamarimhtml/derrida.html} \quad .

[Ze3] Doron Zeilberger, {\it Book Review of J. Borwein and D. Bailey's Mathematics By Experiments and J. Borwein, D. Bailey and R. Girgensohn's Experimentation in Mathematics},
American Scientist {\bf 93 \#2} (March/April 2005), 182-183. \hfill\break
{\tt https://sites.math.rutgers.edu/\~{}zeilberg/mamarim/mamarimPDF/mathexp.pdf} \quad .

[Zu] Wadim Zudilin, {\it Ap\'ery-like difference equations
for Catalan's Constant}, \hfill\break
{\tt https://arxiv.org/abs/math/0201024} \quad .

\bigskip
\hrule
\bigskip
Robert Dougherty-Bliss and Doron Zeilberger, Department of Mathematics, Rutgers University (New Brunswick), Hill Center-Busch Campus, 110 Frelinghuysen Rd., Piscataway, NJ 08854-8019, USA. \hfill\break
Email: {\tt  robert.w.bliss at gmail dot com} \quad, \quad {\tt DoronZeil at gmail dot com}   \quad .

First Written: Sept. 10, 2021. This version: Sept. 21, 2021.

\end